# A finiteness property on monodromies of holomorphic families.

Thomas Delzant

IRMA, Université de Strasbourg et CNRS

## 1. Introduction.

Let $X$ be a Kähler manifold. A holomorphic family of Riemann surfaces of genus $g$ over $X$ is a pair $(Y, \pi)$ where $Y$ is a complex manifold and $\pi\colon Y \to X$ a holomorphic submersion with Riemann surfaces of genus $g$ as fibers. It is called non isotrivial if the family of Riemann surfaces $Y_s = \pi^{-1}(s)$ is not constant in the moduli space of Riemann surfaces. Such a family determines a monodromy, which is a homomorphism $\varphi$ from the fundamental group $\pi_1(X, s_0)$ to the mapping class group $M(S)$ of the topological surface underlying $Y_{s_0}$.

A fundamental result, (and in fact the basic result of the theory), due to Parshin and Arakelov ([Ar] [Par]) and answering a question of Shafarevich asserts that given a Riemann surface $B$ the set of families of given genus over $B$ is finite (for a proof based on the study of Teichmuller spaces, see [Im-Sh], or [McM]). Another consequence of the study of [Ar] is that the number of non isotrivial families over a projective manifold $X$ can be bounded in terms of this manifold.

A uniform result has even been described by L. Caporaso [Ca] who proved that the Hilbert polynomial of a complex surface which is a non singular bundle of genus $p$ over a base of genus $g$ can only take a finite number of values. A consequence is that, given a surface $\Sigma_g$ of genus $g$, and of topological surface $S$ (of genus $p$), the cardinal of the subset the set of homorphisms from the fundamental group of $\Sigma_g$ to the mapping class group of a surface $S$ which can be realized as a monodromy is finite modulo conjugacy at the target and automorphim at the source.

In this article, we want to give a bound in term of the fundamental group $\Gamma$ of the base $X$, independant of the manifold $X$. Before stating our main result let us recall some definitions.

A finitely presented group $\Gamma$ is *Kähler* if it can be realized as the fundamental group of a compact Kähler manifold. The group $\Gamma$ *fibers* if there exists a topological 2-orbifold $\Sigma$ together with a surjective homomorphism $\pi\colon \Gamma \to \pi_1^{\mathrm{orb}}(\Sigma)$ with a finitely generated kernel ; this is equivalent to the fact that every compact Kähler manifold $X$ with fundamental group $\Gamma$ admits a holomorphic map with connected fibers on a complex 1-orbifold with $\Sigma$ as a underlying topological orbifold. Analogously, one says that the family $Y$ factorizes through a curve $B$ if there exists a Riemann orbifold $B$ and a map $q\colon X \to B$ with connected fibers so that $Y$ is the pull-back of some family over $B$. We shall see that this property only depends on the monodromy not on the manifold. The main result of this article is :

**Theorem 1.** *Let $\Gamma$ be a Kähler group, $S$ a topological surface. There exists a finite number of conjugacy classes of homomorphims $\varphi\colon \Gamma \to M(S)$ which can be realized as the monodromy of a holomorphic family of Riemann surfaces on some Kähler manifold with fundamental group $\Gamma$, but do not factorize through a curve.*

Combining this result with the case of curves ([Ca]), one obtain that the number of non isotrivial families over a Kähler manifold $X$ only depends on its fundamental group.

In fact, theorem 1 appears as a (very) special case of a general factorisation theorem (Corollary 11) for actions of Kähler groups on hyperbolic spaces.

Let $H$ be a hyperbolic space in the sense of Gromov ([Gr]), and $G$ a subgroup of the group of isometries of $H$. If $\Gamma$ is a finitely generated group, one can study infinite sequences of non elementary homomorphisms of $\Gamma$ to $G$ with an asymptotic method (sometimes called the Bestvina-Paulin method). Let $\Sigma$ be a fixed set of generator of $\Gamma$ ; the energy of the homomorphism $h$ is $e(h) = \mathrm{Min}_{x \in X} \mathrm{Max}_{s \in \Sigma} d(x, h(s)x)$. An infinite sequence of pairwise non conjugate homomorphisms is diverging if $\limsup e(h_n) = +\infty$. After choosing some ultrafilter, infinite sequences of diverging energy converge to an action of $\Gamma$ on some $\mathbb{R}$-tree : the asymptotic cone of $H$. Due to the fundamental work of Gromov-Shoen, one deduces that the Kähler group "fibers", i.e. admits a epimorphism $\pi$ to a 2-orbifold group, with finitely generated kernel. In many important cases, one further prove that, for an infinity of integers, $h_n$ factorizes through $\pi$.

Unfortunately this simple method cannot be applied directly to the mapping class group acting on the complex of curves (which is a hyperbolic space due to the work of Masur and Minsky [M-M]), as there is no reason that an infinite sequence of homomorphisms to the mapping class group acting on the complex of curves has diverging energy. The whole machinery of [Be-Br-Fu] which creates an action of the mapping class group on a finite product of hyperbolic spaces will be needed, together with the work of [Dr-Be-Sa] on the asymptotic geometry of this group.

In the first paragraph we review the basic results concerning Kähler groups acting on $\mathbb{R}$-trees, in the second paragraph we recall the Bestvina-Paulin method to construct actions on $\mathbb{R}$-trees and we introduce the notion of a limit group associated to an infinite sequence. This theory is applied to the study of Kähler groups. In the last paragraph these results are combined with the description of the asymptotic geometry of the mapping class group to study monodromies.

**Aknowledgements.** I would like to thank Pierre Py for a number of comments and remarks on this subject, and Chloé Perin for discussions concerning limit groups.

## 2. Kähler groups and $\mathbb{R}$-trees.

### 2.1. 2-Orbifold and their fundamental groups.

Let us recall some basic facts about $2-$orbifolds, introduced by W. Thurston ([Th] Chap. 13).

A (oriented) 2-orbifold is a topological surface $S$ (oriented) endowed with a finite set of marked points $\{(p_1, m_1), ...(p_k, m_k)\}$, where $m_i$ is an integer greater than 2. We shall denote it by $\Sigma = \{S; (p_1, m_1), ...(p_k, m_k)\}$. For each point $p_i$, let $\gamma_i$ be a small simple loop which is the boundary of a small embeded disc $D_i$ centered at $p_i$, and choosen so that the disks $(D_i)_{1 \leqslant i \leqslant k}$ are disjoints.

The *fundamental group* of $\Sigma$ is $\pi_1^{\mathrm{orb}} = \pi_1 \left( S \setminus \{p_1, ..., p_k\} \right) / \ll \gamma_1^{m_1}, ..., \gamma_k^{m_k} \gg \}$. The *Euler characteristic* of $\Sigma$ is $\chi(\Sigma) = 2 - 2g - k + \sum_{1 \leqslant i \leqslant k} \frac{1}{m_i}$. A 2-orbifold is called of general type (or hyperbolic) if its Euler characteristic $\chi(X)$ is strictly negative.

If $\Lambda$ is a co-compact lattice in $\mathrm{PSL}(2, \mathbb{R})$, the quotient of the unit disk by $\Lambda$ has a natural structure of a 2-orbifold : $D/\Lambda$ is a topological oriented surface $S$, and modulo $\Lambda$, only a finite number of points $(p_i)_{1 \leqslant i \leqslant k}$ have a non trivial isotropy group, which is a finite cyclic subgroup of order $(m_i)_{1 \leqslant i \leqslant k}$ of $\Lambda$. If $\Sigma = \{S; (p_i, m_i)\}$ is the underlying orbifold, one proves that $\Gamma = \pi_1^{\mathrm{orb}}(\Sigma)$. A *hyperbolic structure* on $\Sigma$ is a realisation of its fundamental group as a co-compact lattice in $\mathrm{PSL}(2, \mathbb{R})$.

A *complex structure* on $\Sigma$ is, by definition, a complex structure on $S$ (marked at the points $p_i$), and the uniformization theorem for 2-orbifolds implies that there is a one-to-one correspondance between hyperbolic and complex structures on general 2-orbifolds.

Let $Y$ be some complex structure on $\Sigma$. In order to define holomorphic maps with values in $Y$, let us choose, for every $i$, a small closed disk $D_i$ around $p_i$. Let $G_i \colon (\widetilde{D}^{m_i}, p_i) \to (D_i, p_i)$ be the $m_i$ – th cover ramified at the origin. If $X$ is a complex manifold, a map $f \colon X \to \Sigma$ is called *holomorphic* if it is holomorphic in the usual sense in $X - f^{-1}\{(p_i)_{1 \leqslant i \leqslant n}\}$, and if for every $i$ and every point such that $f(x) = p_i$, the map $f$ admits a local holomorphic lift through the map $G_i$. In other words, $f$ is locally the $m_i$ – th power of a holomorphic map. This enable one to endow $Y$ with the Kobayashi metric which coincides with its hyperbolic structure.

The following result is probably well known, but does not appear in its full generality in the litterature (see however the article of Catanese [Ca 1], as well as the appendix by Beauville in [Ca2], in the case where the singular set of $\Sigma$ is empty).

**Theorem 2.** *Let $\Gamma$ be the fundamental group of a Kähler manifold $X$, and $\Lambda$ the fundamental group of a general 2-orbifold $\Sigma$. The following are equivalent:*
  i. *There exists a surjective homomorphism $\psi \colon \Gamma \to \Lambda$ with finitely generated kernel.*
  ii. *There exists a complex structure $Y$ on $\Sigma$, together with a holomorphic map $X \to Y$ with connected fibers.*

**Proof.** Let us choose some discrete co-compact action of $\Lambda$ on the unit disk, with orbifold quotient $\Sigma_{\text{aux}}$ (an auxiliary hyperbolic structure on $\Sigma$). The theory of harmonic maps (due to Carlson and Toledo [Ca-To]) proves that there exists an equivariant harmonic map from the universal cover of $X$ to $D$, which descends to a differentiable map $f$ from $X$ to $\Sigma_{\text{aux}}$. Using a Bochner formula, Carleson and Toledo prove that the map $f$ is pluriharmonic, and that the connected components of fibers of $f$ are complex hypersurfaces. The set of connected components of fibers of $f$ is a complex orbifold $Y$ whose singular points are the multiple fibers, with their multiplicity. This orbifold has therefore a continuous map $q \colon Y \to \Sigma$, surjective on fundamental group. In order to prove that this map induces an isomorphism on fundamental groups, we adapt the argument of Catanese. Let us consider the image of $\ker \psi$ in $\pi_1^{\text{orb}}(Y)$. As $\ker \psi$ is f.g. and $\psi' \colon \Gamma \to \pi_1^{\text{orb}}(Y)$ is surjective, the image of $\ker \psi$ in $\pi_1^{\text{orb}}(Y)$ must be trivial or of finite index (finitely generated normal subgroups in 2-orbifold groups are of finite index). This group, which is the kernel of $q_*$, cannot be of finite index as $q_* \colon \pi_1(\Sigma') \to \pi_1(\Sigma)$ is onto, thus it is trivial and $q$ is an isomorphism. $\square$

Using the fact that holomorphic maps are 1-Lipshitz for the Kobayashi metric (which is the hyperbolic metric on a hyperbolic 2-orbifold), one proves (see [Co-Si] for a algebraic proof, [De] Thm.2 for a proof by compacity).

**Proposition 3.** *Let $X$ be a complex manifold. There exists only a finite number of pairs $(Y_i, F_i)$ where $Y_i$ is a hyperbolic/complex 2-orbifold, and $F_i$ a holomorphic map from $X$ to $Y_i$.* $\square$

This suggest the definition (see [ABCKT] Chapter 2, paragraph 3).

**Definition 4.** *A Kähler group $\Gamma$ fibers if there exists a 2-orbifold of general type $\Sigma$ and a surjective homomorphism with finitely generated kernel $\Gamma \to \pi_1(\Sigma)$.*

We wish to emphazise that this definition only depends on $\Gamma$, not on the choice of a Kähler manifold with $\Gamma$ as a fundamental group. Combining the two previous theorems, one gets the following :

**Proposition 5.** *Let $\Gamma$ be a Kähler group. There exist a finite family of pairs $(\Sigma_i, \pi_i)_{1 \leqslant i \leqslant p}$ of 2-orbifolds of general type and surjective homomorphisms $\pi_i \colon \Gamma \to \pi_1^{\text{orb}}(\Sigma_i)$ with finitely generated kernel, such that for every Kähler manifold $X$ with fundamental group $\Gamma$, every hyperbolic orbifold $Y$ and every holomorphic map $F \colon X \to Y$ with connected fibers, there exists and integers $i$ such that the underlying orbifold of $Y$ is isomorphic to $\Sigma_i$ by an isomorphism such that the induced map on fundamental group is $\pi_i$.* $\square$

**2.2. Trees.**

Recall that an $\mathbb{R}-$tree is a connected geodesic metric space which is $0-$hyperbolic (see [Be1] for an introduction to this subject).

The (isometric) action of a group on an $\mathbb{R}-$tree is called *minimal* if it does not contain an invariant sub-tree. If the group is finitely generated, it admits a unique minimal invariant sub-tree. The action is called *non-elementary* if it is neither *elliptic* (fixing a point) nor *axial* (fixing a line but no point on this line) nor *parabolic* (fixing a unique point at infinity).

The main example of a Kähler group acting on a $\mathbb{R}-$tree is the fundamental group of a Riemann surface (or 2-orbifold) endowed with a holomorphic quadratic differential $\omega$ : the $\mathbb{R}-$tree is the set of leaves of the real part of $\omega$. (The notion of a quadratic differential can be extended to the case of an orbifold : let $S$ be the underlying Rieman surface, $K$ the canonical divisor and $D$ the divisor $D = \Sigma(m_i - 1)p_i$. Then a quadratic differential is a holomorphic section of $2K + 2D$.

After the fundamental work of Gromov-Schoen [G-S] several authors ([K-S], [Su]) studied actions of Kähler groups on $\mathbb{R}-$trees, and the following theorem summarizes the situation.

**Theorem 6.** *Let $\Gamma$ be a Kähler group acting on a $\mathbb{R}$-tree $T$. Assume that $T$ is minimal and is not a line. Then $\Gamma$ fibers. Moreover, there exists a 2-orbifold $\Sigma$ of general type, a surjective map $\pi \colon \Gamma \to \pi_1^{\mathrm{orb}}(\Sigma)$ with finitely generated kernel, an action of the fundamental group of $\Sigma$ on an $\mathbb{R}-$tree $T'$ and an $\pi-$equivariant map $T' \to T$. If the action of $\Gamma$ on $T$ is faithfull, then $\pi$ is an isomorphism.*

If $T$ is a simplicial locally finite tree, the proof is explained in [G-S] ; the general case is sketched in the paragraph 9.1 of the same article. Let us recall the main steps of the proof.

Let $X$ be a Kähler manifold with fundamental group $\Gamma$. One construct a $\Gamma$-equivariant harmonic map of finite $\Gamma-$energy $h$ from the univeral cover of $X$ to the tree $T$ ([G-S], [K-S]). The regularity of this harmonic map ([G-S]) has been detailed by [Su] ; in particular the set of singular points $Y$ is of codimension 2, and the image of a connected fundamental domain in $\tilde{X}$ appears to be a finite tree (the convex hull of a finite number of points). At this point one can copy the argument of [G-S] : using the Kähler structure, one prove that the map $h$ is pluriharmonic : outside from $Y$ it is locally the real part of a holomorphic function, and there exists a holomorphic quadratic differential $\omega$ on $X$ such that locally $\pm dh = \omega$. Using the fact that the tree is not a line one proves that one leaf of the foliation defined by $\omega$ on $X$ is singular (entirely contained in the set $\omega = 0$) hence compact. One deduces that all leaves are compact and $X$ fibers on some hyperbolic 2-orbifold $\Sigma$ with fibers the leaves of $\omega$. Thus $\omega$ comes from $\Sigma$, and if $T_\omega$ denote the $\mathbb{R}-$tree dual to the leaves of $\omega$, the harmonic map $h$ factorises through the $\pi_1(\Sigma)$ equivariant map from $D \to T_w$ . $\square$

The case of axial actions of Kähler groups is also quite usefull (see [De2]).

**Theorem 7.** *Let $\Gamma$ be a Kähler group with an isometric action on a line (isometric to $\mathbb{R}$). If the kernel of this action is not finitely generated, the group $\Gamma$ fibers ; more precisely, there exists a 2-orbifold of general type a surjective morphism with finitely generated kernel $\pi \colon \Gamma \to \pi_1^{\mathrm{orb}}(\Sigma)$ such that the action factorizes through $\pi$.*

Let $X$ be a compact manifold with fundamental group $\Gamma$. In the case of an oriented line, such an action is defined by a holomorphic form, and [De2] applies. In the general case, the action defines a holomorphic quadratic differential $\omega$ which becomes a differential on a (perhaps ramified) double cover. All the leaves of this form are compact, hence all the leaves of $\omega$ are compact and the same argument as in Theorem 5 concludes. $\square$

**3. Limit groups and $\mathbb{R}$-trees as limit spaces.**

In a number of cases, actions of group on ℝ-trees are obtained by a limit process of actions on hyperbolic spaces, as first observed by Bestvina [Be1] and Paulin [Pa].

### 3.1. Limit spaces and limit groups.

Let $\Gamma$ be a finitely generated group, $(H_n, x_n, d_n)$ a sequence of pointed metric spaces, $\varphi_n$: $\Gamma \to \mathrm{Isom}(H_n)$ a sequence of isometric actions of $\Gamma$ on $H_n$.

Let us fix a non principal ultrafilter $\omega$ on $\mathbb{N}$. Recall that $\omega$ is a subset of the set of subset of $\mathbb{N}$ such that $\varnothing \notin \omega$, if $A, B \in \omega, A \cap B \in \omega$, the complementary $F^c$ of every finite set is in $\omega$, and $\omega$ is maximal for these properties. Recall (Bolzano-Weierstrass) that every bounded sequence has a $\omega-\mathrm{limit}$ : for every $\varepsilon > 0$ the set $\{n \in \mathbb{N} / |u_n - l| < \varepsilon\}$ belongs to $\omega$.

**Definitions.** The limit space $H_\omega$ is the quotient of the set $\Pi^{\mathrm{bounded}} H_n = \{(y_n) / \lim_\omega d_n(y_n, x_n) < \infty\}$ by the equivalence relation $\lim_\omega d_n(y_n, z_n) = 0$. This space has a natural base-point, the equivalence class of $(x_n)$, and a natural distance $d_\omega$ defined by $d_\omega(y,z) = \lim_\omega d(y_n, z_n)$.

If for every element in a set of generators of $\Gamma$, the sequence $d_n(x_n, gx_n)$ is bounded, the group $\Gamma$ acts isometrically on $X_\omega$.

The *stable kernel* of $(\varphi_n)_{n \in \mathbb{N}}$ is $K_\omega = \{g \in \Gamma / \{n / \varphi_n(g) = e\} \in \omega\}$, and the *limit group* $\Gamma_\omega$ is the quotient $\Gamma / K_\omega$.

We note that, by definition, $\Gamma_\omega$ acts on $H_\omega$. Furthermore, the following proposition will be usefull. It emphazise the importance of finitess properties of kernels of actions.

**Proposition 8.** *Let $K$ be a finitely generated subgroup contained in $K_\omega$. Then for some $n$ (and in fact for $\omega$ almost all $n$), $K \subset \mathrm{Ker}(\varphi_n)$ the kernel of the action $\varphi_n$.* □

### 3.2. Hyperbolic metric spaces and trees.

Let $H$ denote a $\delta$-hyperbolic space (in the sense of Gromov [Gr]) , $\Gamma$ a finitely generated group, $\sigma$ be a fixed set of generators of $\Gamma$, and $\varphi: \Gamma \to \mathrm{Isom}(H)$ an isometric action.

The *energy* of a point $x$ is defined by $e(x) = \mathrm{Sup}_{s \in \sigma} d(sx; x)$, and the *energy* of the action is the minimum $e(\varphi) = \mathrm{Min}\ e(x)$. If the action of $\Gamma$ on $H$ is not elementary (neither elliptic, nor parabolic, nor loxodromic), then the fonction $e(x)$ is metrically proper, and the set $\{x / e(x) \leqslant e(h) + 1\}$ is not empty of bounded diameter. If $H$ is CAT($-1$) complete, the minimum is in fact achieved, as $e$ is a convex proper function. (See [Be2]).

The main result concerning limits of actions of a group acting on a hyperbolic space is due to Bestvina and Paulin see ([Be1], [Pa], [Be 2], [K-S]).

**Theorem 9.** *Let $\Gamma$ be a finitely generated group, $(H, d)$ a hyperbolic space and $\varphi_n$: $\Gamma \to$ Isom $(H)$ be a sequence of actions whose energy $e_n$ goes to infinity. Let $x_n$ be choosen so that $e_n(x_n) \leqslant e_n + 1$, and $d_n = \frac{d}{e_n}$. Then the limit space $\lim_\omega (H, d_n, x_n)$ is a complete $\mathbb{R}-$tree $H_\omega$, the action of $\Gamma$ on $H_\omega$ has energy 1, and the minimum of energy is achieved at the origin $x_\omega$.*

Proof : see [Be2], [K-S]. □

Remark. Sometimes, the space $H$ is a smooth (perhaps infinite dimensional) manifold. In this case, the work of [K-S] directly produces an harmonic map with values in the limit space, as a limit of harmonic maps on balls in $(H, d_n, x_n)$ centered at $x_n$ and with radius $r_n \to \infty$.

As the energy is 1, this action cannot be elliptic ; it therefore contains a minimal invariant subtree which will be denoted by $T_\omega \subset H_\omega$. In many cases, one knows something about the kernel of the action of $\Gamma_\omega$ on $T_\omega$. The following definitions are useful in our context.

*Acylindricity.* One says that the action of $G$ on a $\delta$-hyperbolic space $H$ is $K$-acylindrical if there exists an integer $N$ such that for every pairs of point $a, b$ with $d(a, b) \geqslant K$ the set $\{g \in G /\ d(ga, a) + d(gb, b) \leqslant 1000\delta\}$ is finite of cardinal $\leqslant N$.

Let the group $G$ acts on a $\delta$-hyperbolic space $H$. The translation length $[h]$ of an element $h$ is $[h] = \text{Min}_{x \in H} d(x, hx)$. Let $h$ be some element with translation length $\geqslant 100\delta$. Then $h$ is hyperbolic : if $x_0$ realizes the minimum $\text{Min}_{x \in H} d(x, hx)$, the sequence $(h^n . x_0)_{n \in \mathbb{Z}}$ is a $100\delta$ local geodesic) (see for instance [CDP], chap. 9 for $6\delta$ instead of $100\delta$). Let $h^\pm$ the two fixed points of $h$ at infinity, $L_h = \{x / d(x, hx) \leqslant [h] + 10\delta\}$. This set $L_h$ is an invariant tubular neigbourhood the union of lines joining the two fiwed points of $h$ at infinity. Let $E(h)$ be the elementary group generated by $h$ (the set of $g$ such that $g\{h^+, h^-\} = \{h^+, h^-\}$).

*Weak Weak Proper Discontinuity.* Definition ([BBF1,2]) On says that the action of $G$ on a hyperbolic graph is *WWPD* if, for every hyperbolic element $g$, there exists a $D$ and a subgroup $C(g) \subset E(g)$ such that for every $h \notin C(g)$ the diameter of the projection of $h.L_g$ on $L_g$ is bounded by $D$.

**Theorem 10.** *Let $(H, d)$ be a hyperbolic space, $G$ a group of isometries of $H$, and $h_n \colon \Gamma \to G$ a sequence of actions whose energy $e_n$ goes to infinity. If the action of $G$ on $X$ is acylindrical, the kernel of the action of $\Gamma_\omega$ on $T_\omega$ is finite if $T_\omega$ is not a line, or virtually abelian if $T_\omega$ is a line.*

*If the action is WWPD, and $K \subset G$ is a f.g. subgroup whose image is in the kernel of the action of $\Gamma_\omega$ on $T_\omega$, then for an infinity of $n$, $h_n(K)$ is contained in $E(h_n(h))$, the maximal elementary subgroup containing some hyperbolic element $h_n(h)$ of $G$ ; if $T_\omega$ is not a line, for an infinity of $n$, the group $h_n(K)$ is elliptic.*

Proof. As $\Gamma_\omega$ is finitely generated, and the action of $\Gamma_\omega$ is not elliptic, one can find an element $h$ which acts as a hyperbolic isometry on $T_\omega$ : $T_\omega$ contains a $h$-invariant line $\Lambda_h$. Let $K$ be a finitely generated subgroup of $\Gamma$ which fixes this line, and $k \in K$. Using the WWPD property, one gets that, for $\omega$-almost integer $n$, the subgroup generated by commutators $h_n([h, k])$ is contained in $C(h_n(h)) \subset E(h_n(h))$. Therefore $h_n(khk^{-1}) \in E(h_n(h))$, and $h_n(k) \in E(h_n(h))$. Thus, as $K$ is finitely generated, $\omega$ − almost surely, $h_n(K) \subset E(h_n(g))$.

If the action is acylindrical, the subgroup generated by commutators $h_n([k, h])$ is finite for $n \gg 1$, and the kernel of the action $h_n$ is therefore virtually abelian.

If $T_\omega$ is not a line, for some $u$, $\Lambda_{uhu^{-1}} = h_\omega(u)\Lambda_g \neq \Lambda_g$, and for the same reason, if $H$ is a subgroup of the kernel of the action, $h_n(H) \subset E(h_n(uhu^{-1}))$. The group $h_n(H)$ is contained in two elementary subgroups generated by $h_n(h)$ and $h_n(uhu^{-1})$ is elliptic. And therefore finite in the acylindrical case. □

### 3.3. Kähler groups.

Let us keep the notations of 2.2 unchanged, but assume now that $\Gamma$ is a Kähler group. Applying Theorems 10 and 6 (or 7 for the axial case) simultaneously one gets.

**Corollary 11.** *Let $\Gamma$ be a Kähler group; let $G$ a group of isometries of a hyperbolic space $(H, d)$, $h_n \colon \Gamma \to G$ a sequence of actions whose energy $e_n$ tends to infinity and $x_n$ be chosen so that $e_n(x_n) \leqslant e_n + 1$. Let $X_\omega = \lim_\omega (\frac{1}{e_n} X_n, x_n)$, and $T_\omega$ be the minimal invariant subtree in $X_\omega$.*

*If $T_\omega$ is not a line, then $\Gamma$ fibers. There exists a 2-orbifold $\Sigma$ and a surjective homomorphism $\pi \colon \Gamma \to \pi_1^{\text{orb}}(\Sigma)$ with finitely generated kernel $N$. The same is true if $T_\omega$ is a line, but the kernel of the action of $\Gamma$ on this line is not finitely generated.*

*If furthermore, the action of $G$ on $H$ is acylindrical, then for an infinity of n, $h_n(\ker \pi)$ is finite group if $T_\omega$ is not a line or is a line but the kernel is infinitely generated, and $\pi$ factorizes through a finite extension of $\pi_1^{\mathrm{orb}}(\Sigma)$. It factorise through a virtually abelian group if $T_\omega$ is a line and the kernel is finitely generated.*

*If the action is WWPD, then for an infinity of n, the group $h_n(N)$ is contained in the elementary subgroup of a hyperbolic element of $G$.*

## 4. Mapping class group, complex of curves and holomorphic families of Rieman surfaces.

### 4.1. Asymptotic geometry of the mapping class group.

Let $S$ be a topological surface (maybe with a boundary) ; the complex of curve (and arcs) of $S$, denoted by $X(S)$, or $X$ if only one surface is involved, is the graph whose vertices are the homotopy classes of simple closed loops on $S$ non parallel to the boundary, and arcs from the boundary to the boundary but not homotopic to the boundary, and whose edges are pair of curves who can be made disjoint by a homotopy.

A fundamental result, due to H. Masur and Y. Minsky asserts that $X$ is hyperbolic [M-M] ; this statement has been improved by B. Bowditch [Bo] who proved that the action of $M(S)$ on $X$ is acylindrical. Therefore, in principle, the results in paragraph 2 can be applied to this example. However, in the study of monodromy groups we will need a deeper result, and study sequences in Hom $(\Gamma, M(S))/$ conj which are infinite but not of diverging energy when acting in $X(S)$. The recent work of M. Bestvina, M. Bromberg, and K. Fujiwara [Be-Br-Fu 1,2] (introducing the WWPD property) gives a new geometrisation of the mapping class group which enables to apply the results of paragraph 3, and in particular the corollary 11.

**Theorem 12.** *[Be-Br-Fu 1,2] The mapping class group $M(S)$ of a surface $S$ contains a subgroup of finite index $M_1(S)$ which admits a product action on a finite product of hyperbolic spaces $X = \Pi_{i \in I} X_i$, each of these action being faithfull and WWPD. Furthermore the orbit map $M(S) \to X$ is a quasi-isometric embedding.*

The space $X$ is not hyperbolic, but any of its asymptotic cones is a finite product of $\mathbb{R}$-trees. We want to explain how to modify the argument of 3.2 to create actions on $\mathbb{R}-$ trees from an infinite sequence of pairwise non conjugate homomorphisms of a finitely generated group $(\Gamma, \Sigma)$ to $M(S)$. At this point one could directly apply the work of J. Behrstock, C. Drutu and M. Sapir [Be-Dr-Sa] Cor.6 of the appendix, but as we need to get more informations on the kernel of the action on the limit tree, we prefer to deduce our result from [Be-Br-Fu].

In order to apply corollary 11 to this action, we recall a fundamental result of Ivanov [Iv] concerning the mapping class group : *a finitely generated subgroup of $M(S)$ which is not reducible must contain a pseudo-Anosov element.* The following lemma follows.

**Lemma 13.** *Let $G \subset M(S)$ be a f.g. irreducible group, $N \subset G$ an infinite f.g. normal subgroup. Then $N$ must contain a pseudo-Anosov.*

Proof. If $N$ contains no-pseudo Anosov, it is elliptic acting in the complex of curves, and the set of curves preserved by $N$ is bounded. As $N$ is normal, this set is preserved by $G$ which is therefore elliptic acting on this complex. Applying Ivanov's theorem, we get that $G$ must be reducible.□

Now, let $\Gamma$ be a group and $h_n\colon \Gamma \to M(S)$ be an infinite sequence of pairwise non conjugate homomorphisms. As $I$ is finite, $\Gamma$ admits a subgroup of finite index $\Gamma_1$ such that the restriction of $h_n$ to $\Gamma_1$ do not permute the factors $X_i$. Let us fix a generating system $\Sigma$ of $\Gamma$. Simultaneously, we fix a generating system of $M(S)$ and denote $|\gamma|$ the word length of an element $\gamma$.

Let $\alpha_0 \in X$ be a fixed base-point in $X$. It is know ([Be-Br-Fu]) that the orbit map $\psi\colon M(S) \to X$ defined by $\psi(\gamma) = \gamma.\alpha_0$ is a quasi-isometric embedding. In other words, there exists a constant $K$ such that : $\quad K^{-1}|\gamma| - K \leqslant d(\alpha_0, \gamma.\alpha_0) \leqslant K|\gamma|$.

On the other hand, the family $h_n$ is pairwise non conjugate, therefore the sequence $E_n = \mathrm{Min}_{\gamma \in M(S)} \mathrm{Max}_{s \in \Sigma} |\gamma^{-1} h_n(s) \gamma|$ is unbounded, as well as the sequence $e_n = \mathrm{Min}_{\gamma \in M(S)} \mathrm{Max}_s d(h_n(s).\gamma\alpha_0, \gamma\alpha_0)$ as $K^{-1}E_n - K \leqslant e_n \leqslant KE_n$. After conjugating $h_n$, we may assume that this Min is achieved for $\gamma = 1$.

Let $C = \lim_\omega (\frac{1}{E_n} \mathrm{Ca}(M(S), e))$ be the asymptotic cone of the mapping class group associated to this divergent sequence. The orbit map (for some fixed origin in $\Pi_{i \in I} X_i$) induces an equivariant bilipshitz embedding of $C$ in a product of tree $\Pi_{i \in I} T_i$, with $T_i = \lim_\omega \frac{1}{e_i} X_i$.

**Proposition 14.** *Assume that for all $n$, the group $h_n(\Gamma)$ is neither reducible, nor virtually abelian. Let $N$ be a f.g. normal subgroup contained in the kernel of the action of $\Gamma$ on $T_i$ ; then for an infinity of $n$, $h_n(N)$ is finite, and $\Gamma$ has a subgroup of finite index $\Gamma_1$ such that the restriction of $h_n$ to $\Gamma_1$ factorizes through $\Gamma_1/N$.*

Proof. The first fact (for some $i$, the action of $\Gamma$ on $T_i$ is not ellliptic) directly follows from [Be-Dr-Sa], Section 6, Thm. 6.2, who proved that the orbit of the limit group $\Gamma_\omega$ on $C$ is unbounded. As the embedding of $C$ in this product is bilipshitz (thus metricaly proper) it cannot project onto bounded orbits. As $\Gamma$ is finitely generated, we can find an element $h \in \Gamma_1$ and an index $i$ such $\lim_\omega (h_n(h))$ is hyperbolic, and acting in $X_i$, $h_n(h)$ is hyperbolic for almost all $n$. Due to the WWPD property, $h_n(N) \subset E(h_n(h))$. If $h_n(h)$ is reducible $E(h_n(h))$ is reducible too. Applying Lemma 13, we get a contradiction. Therefore $h_n(h)$ is a pseudo-Anosov, $E(h_n(h))$ is virtually $\mathbb{Z}$, and cannot contain a infinite normal subgroup of $h_n(G)$, unless $h_n(G) \subset E(h_n(h))$ is virutally abelian. Thus $h_n(N)$ is finite. $\square$

### 4.2. Holomorphic family of Riemann surfaces.

Let $X$ be a compact Kähler manifold. A holomorphic family of Riemann surfaces of genus $g$ is a pair $(Z, \pi)$ where $Z$ is a compact complex surface, and $\pi\colon Z \to X$ a holomorphic fibration such that the fibers of $\pi$ are Riemann surfaces of genus $g$.

Let $x_0$ be a base point in $X$, $Z_{x_0} = \pi^{-1}(x_0)$ its fiber, and $S$ be the underlying topological surface (with punctures). Let $M_{g,p} = M(S)$ be the mapping class group of $S$, and $\mathcal{T} = \mathcal{T}_{g,p}$ the corresponding Teichmuller space. Note that $Z_{x_0}$ determines a natural base point $y_0$ in $\mathcal{T}$, and an action of $M(S)$ on $\mathcal{T}$.

These notations being fixed, the family $\pi$ determines two objects.

a. the holonomy : $\varphi\colon \pi_1(X, x_0) \to M(S)$
b. the classifying map $\Phi\colon (\tilde{X}, x_0) \to (\mathcal{T}, y_0)$ from the universel cover of $X$ to the Teichmuller space, which is a $h$-equivariant holomorphic map.

One says that the family is isotrivial if $\Phi$ is constant (equivalently, the image of $\varphi$ is finite, as the action of $M(S)$ on $\mathcal{T}$ is proper with finite stabilizers).

The existence of a "Universal Teichmuller curve" shows that conversely a pair of such maps determines a family over $X$ (in the case of finite holonomy, one uses the solution by Kerchoff to the Nielsen realization problem).

It is useful, for our purpose, to extend the definition to the case where the base is a complex (or hyperbolic) 2-orbifold. Let $\Gamma \subset \mathrm{PSL}(2, \mathbb{R})$ be a fuchsian group, and $X = U/\Gamma$ the correponding complex 2-orbifold. A family of Riemann surfaces overs $X$ is given by a representation $h: \Gamma \to M$ and a holomorphic $\Gamma$-equivariant map from $U \to \mathcal{T}$. A priori this definition can also be given for a base which is a projective line or an elliptic curve ; but as the Teichmuller space is hyperbolic in the sense of Kobayashi, the classifying map over such a base is always constant and the family isotrivial.

**Definition 15.** *Let $(Z, \pi, X)$ be a family of Riemann surfaces over X. One says that $\pi$ factorizes through a curve if there exits a hyperbolic 2-orbifold $Y$ a holomorphic map with connected fibers $F: X \to Y$ and a family $(Z', \pi', Y)$ over Y such that $Z = F^* Z'$.*

Let us recall that the Teichmuller space is contractible, and admits an invariant Kähler form $\Omega_{\mathrm{WP}}$, the imaginary part of the Weil-Peterson Kähler metric. As $\mathcal{T}$ is contractible, and the action of $M$ is proper with finite stabilizers, the form $\omega$ determines a cohomology class $[\Omega_{\mathrm{WP}}] \in H^*(M(S), \mathbb{R})$, the Weil-Peterson class.

Keeping these notations we have:

**Proposition 16.** *The following proposition are equivalent :*
  *i. The family factorizes through a curve*
  *ii. $\varphi^*([\Omega_{\mathrm{WP}}])^2 = 0$ ($\in H^4(\Gamma, \mathbb{R})$)*
  *iii. The complex rank of the (holomorphic) classifiying map $\Phi$ is 1.*

Proof. $i \Rightarrow ii$ and $iii \Rightarrow i$ are obvious. Assume that the complex rank of $\Phi$ is $r \geqslant 2$. The image of $X$ in $M_g$ is a compact analytic manifold $P$ on which $\Omega^r$ induces a volume form outside a set of codimension $> 1$. Thus $\int_P \Omega_{\mathrm{WP}}^r \neq 0$ and $\Phi^*[\Omega_{\mathrm{WP}}^r] \in H^r(X, \mathbb{R}) \neq 0$. But, as the Teichmüller space is contracible, the map $\Phi^*: H^*(M, \mathbb{R}) \to H^*(X, \mathbb{R})$ factorizes through $H^*(\Gamma, \mathbb{R})$. □

**Corollary 17.** *The fact that a family of Riemann surfaces factorizes only depends on its monodromy (the morphism $\varphi: \Gamma \to M$), not on the variety.*
*If there exists a finite cover $X_1$ of $X$ such that the pull back of these bundle $Z$ factorize through a curve, then $Z$ itself factorize.*

The following theorem, due to Imayoshi and Shiga [Im-Sh] (see also [McM]) is important.

**Theorem 18.** *Let $h$ be the monodromy of a family of Riemann surfaces ; then the image of $h$ cannot be reducible or virtually cyclic .*

The first point is just a reformulation of the paragraph 4, case 2 in [Im-Sh], or the part "Irreducibility" in the proof of [McM]. In order to check that the image cannot be $\mathbb{Z}$, note that, by contradiction, this would implies that $\varphi^*\Omega = 0$, as $H^2(\mathbb{Z}, \mathbb{R}) = 0$, and $\Phi$ would be constant. □

### 4.3. Finiteness of monodromies.

Let us say that a morphism $\Gamma \to M$ from a Kähler group to the mapping class group is a monodromy if it can be realized as the monodromy of a family of Riemann surfaces over Kähler manifold with fundamental group $\Gamma$.

**Theorem 19.** *Let $\Gamma$ be a Kähler group. Then there exists only a finite number of conjucacy classes of monodromy $\varphi \colon \Gamma \to M$ which do not factorise through a Riemann surface.*

Proof. Let $\varphi_n$ be an infinite sequence of pairwise non conjugate monodromies. Applying the results of paragraph 4 and 3, one construct a subgroup of finite index of $\Gamma$ which fibers over a Riemann surface group $S$, such that if $N$ is the kernel of this fibration, $\varphi_n|_N$ is finite. Thus, $\Gamma$ admits a subgroup of finite index such that $\varphi_n$ restricted to this subgroup factorizes. By Cor. 17, $\varphi_n$ itself factorizes through a Riemann surface, contradiction. □

**Remark 20.** If $X$ is a projective manifold, it contains a Riemann surface $L$ such that the fundamental group of $L$ surjects onto that of $X$. The monodromy of a family over $X$ is well defined by the monodromy of the restriction of this family to $L$, hence can only take a finite number of values. However the bound we get by this argument depends on $X$, even if the fundamental group is fixed.

## 5. Bibliography.

Thomas Delzant 11

IRMA, Université de Strasbourg et CNRS,
7 rue René Descartes
67084 Strasbourg Cedex France
delzant@unistra.fr